\documentclass[14pt]{article}
\usepackage[english]{babel}
\usepackage{amsfonts,amssymb,amsmath,amstext,amsbsy,amsopn,amscd,amsthm,graphicx,euscript,graphicx}
\usepackage{graphics}
\newtheorem{thm}{Theorem}
\newtheorem{lm}{Lemma}

\newcommand{\eps}{\varepsilon}

\newcounter{tdfn}
\newenvironment{dfn}
{\vspace{0.15cm}{\bf Definition \arabic{tdfn}.}} {\par
\addtocounter{tdfn}{1}}
\newcounter{trk}
{\endtrivlist}

\def\:{\colon}

\def\0{{\mathbf 0}}
\def\1{{\mathbf 1}}

\def\C{{\mathbb C}}

 \author{V.O.Manturov \footnote{Moscow Institute of Physics and Technology}}

\title{Realisability of $G_{n}^{3}$, realisability projection, and kernel of the $G_{n}^{3}$-braid presentation}

\begin{document}

\maketitle

\begin{abstract}
The aim of this article is to prove that the kernel of the map from the pure braid group $PB_{n},n\ge 4$ to the
group $G_{n}^{3}$ consists of full twist braids and their exponents.

The proof consists of two parts. The first part which deals with $n=4$ relies on
the crucial tool in this construction having its own interest is the {\em realisability projection} saying
that if two {\em realisable} $G_{4}^{3}$-elements are equivalent then they are equivalent by a sequence of realisable ones.

The second part (an easy one) uses induction on $n$.

\end{abstract}

AMS MSC: 20F36, 57K10, 57M25, 57M27.

\section{Introduction}

In \cite{Gnk}, the author introduced a series of groups $G_{n}^{k}$ for $n>k$ and formulated the following principle:

{\em If dynamical systems describe a motion of $n$ particles and there exists a good codimension $1$ property governed by exactly $k$ particles
then these dynamical systems admit an invariant valued in the group $G_{n}^{k}$}.

The groups are defined by the following presentation:

\begin{dfn}\label{def:gnk}
The {\em $k$-free braid group} $G_{n}^{k}$ is defined by generators
$a_{m}$, $m\subset\{1,\dots,n\}, Card(m)=k$  (so that there are totally $n\choose k$ generators),
and relations:
\begin{enumerate}
\item $a_m^2=1$;
\item (far commutativity): $a_{m}a_{m'}=a_{m'}a_{m}$ for any $m,m'$ with $Card(m\cap m')<k-1$ ;
\item (tetrahedron relation
also known as Zamolodchikov relation,
higher version of the Yang-Baxter equation):
$a_{m^1}\cdot a_{m^2}\cdots a_{m^{k+1}}= a_{m^{k+1}}\cdots a_{m^2}\cdot a_{m^1}$ for any $(k+1)$-tuple $U$ of distinct indices $u_{1},\dots, u_{k+1} \in \{1,\dots, n\}$ and $m^{j}=U\setminus\{u_{j}\}, j=1,\dots, k+1$.
\end{enumerate}
\end{dfn}

Note that the relations of the first type makes our group similar to Coxeter groups \cite{Cox}; the relations of the second type makes it similar
to braid groups \cite{book}. Note that the relations of the second type are void for $n=k+1$.

The simplest example mentioned in \cite{Gnk} and calculated explicitly in \cite{MN}, dealt with (pure) braids
(dynamical systems of $n$ points on the plane) and the following property for $k=3$: three points are collinear
(we require the genericity condition for these collinear triples).

The realisation of this principle led to the homomorphism $f:PB_{n}\to G_{n}^{3}$ which is constructed as follows.

We consider a braid $\beta$ as a motion of $n$ points $z_{1}(t),\cdots, z_{n}(t),t\in [0,1]$.
We assume a pure braid $\beta$ to be {\em generic} (for more details see \cite{MN}) with the initial
points and final points uniformly distributed on the unit circle: $z_{j}(0)=z_{j}(1)=exp(\frac{2\pi j\cdot i}{n}),j=1,\cdots, n$.
We call this collection of points $Z_{0}= \{z_{j}(0)\}$ {\em the  initial configuration}.

As $t$ changes from $0$ to $1$, there are some moments for which some three points $z_{p}(t),z_{q}(t),z_{r}(t)$ are collinear.
For each such moment we write down the generator $a_{pqr}$ (note that by definition $a_{pqr}=a_{qrp}=a_{qpr}$). The product of such generators is our $f(\beta)$.

One of the goals of the present paper is the following
\begin{thm}
The kernel of the map $f:PB_{n}\to G_{n}^{3}$ consists of ``full twists''. \label{maintheorem}
\end{thm}

Here by a full twist we mean the pure braid $Z(t)=\{z_{j}(t)\},j=1,\cdots, n$, where $z_{i}(t) = exp(2\pi \cdot i \cdot t)\cdot z_{i}(0)$.
Obviously, since all points always belong to the same circle, there is no moment when three points are
collinear.

The main tool in proving this theorem relies on the following techniques: we distinguish
between (algebraically) {\em realisable} $G_{n}^{3}$-words (more or less those coming from braids) and
non-realisable ones, see definition ahead.

Our strategy will be: {\em restore a braid from a word} for $G_{n}^{4}$ and then use induction to prove theorem for $G_{n}^{k}$.

\subsection{Acknowledgements}

I am extremely grateful to I.M.Nikonov for pointing out some mistakes in previous versions and to
Huyue Yan for discussing the text and pointing to some remarks.

\section{Realisability}

Having a generic dynamics of points $\beta=\{z_{i}(t)\}$ (with initial and final positions
uniformly distributed on the unit circle: $z_{k}(0)=z_{k}(1)=exp(\frac{2\pi k\cdot i}{n}))$,
we can write down the corresponding word $w=w(\beta)$.

Equivalent dynamics (isotopic braids) give rise to words equal in $G_{n}^{3}$.

However, it is almost obvious that not all words originate from braids. For example, assume that $k=4$ and
for some moment $t$ we have $z_{4}(t)=0, z_{j}(t)=exp(\frac{2\pi j \cdot i}{3}),i=1,2,3$. From this position,
there is no way to move $z_{1}(t),z_{2}(t),z_{3}(t)$ to the collinear position
because there is a point $z_{4}$ inside the triangle. We shall analyse this in more details in later sections.

Denote the set of words in $G_{n}^{3}$-generators by $W$.

We are going to define two types of realisability, {\em algebraic realisability} and {\em realisability by lines}.

We say that a word $w\in W$ is {\em realisable by lines} if $w=w(\beta)$ for some braid $\beta$.

A much weaker notion is {\em algebraic realisability}. Take a word $w\in W$ and let us read it from the left to the
right.

Now, consider a word $w\in W$. Let us define the {\em initial indices} to be the indices
for the points $exp(\frac{2\pi j\cdot i}{N})$, namely, we define
$(i,j,k)=1$ if and only if $mod(j-i,n)<mod(k-i,n)$.

Starting from these fixed indices, we shall transform them as follows:
we read the word $w$ from the left to the right, and whenever we meet the
generator $a_{m}$, $m=\{i,j,k\}$ in whatever order then we change $(i,j,k)$ to
its negative (all permutations of $(i,j,k)$ are changed accordingly).

It is quite easy to see that the above procedure defines an {\em action} of the group $G_{n}^{3}$
on the set of indices: indeed, for each $G_{n}^{3}$-relation the LHS changes exactly the same indices
as the RHS.

Meeting a generator $a_{ijk}$ we say that this letter $a_{ijk}$ in the word $w$ is realisable for the order
$i_{1}<i_{2}<i_{3}, \{i_{1},i_{2},i_{3}\}=\{i,j,k\}$ if $\forall p\neq i,j,k: (i,j,p)=(i,k,p)=(j,k,p)$.
The generator $a_{ijk}$ is realisable if it is realisable for some ordering. Note that
realisability with respect to $i<j<k$ is the same as realisability with respect to $k<j<i$,
so it is reasonable to talk about {\em realisability with respect to the central letter} (say, $j$).

Note that in the case $n=4$ the realisability condition $(i,j,p)=(i,k,p)=(j,k,p)$
has to be checked for exactly one $p$, namely, $p=\{1,2,3,4\}\backslash \{i,j,k\}$.

We say that the word $w$ is {\em (algebraically) realisable} if all letters in it are algebraically realisable
(each with respect to some central letter).

It is not hard to see that a word is {\em realisable by lines} then it is {\em algebraically realisable}.

Indeed, having three points $z_{i},z_{j},z_{k}$ almost on the same line with $i$ in the centre
we see that $\forall p\neq i,j,k: (i,j,p)=(i,k,p)=(j,k,p)$.

The reader can easily find an example when the inverse is not the case.

Later on, by {\em realisability} of  we mean exactly {\em algebraic realisability}, which will be sufficient for our purposes.

\subsection{The key theorem}

The theorem playing a crucial role in the whole construction is the following:

\begin{thm}
Let $w,w'$ be two realisable words representing elements in $G_{4}^{3}$. Then there is a
sequence $w=w_{0}\to \cdots \to w_{k}=w'$ such that any two adjacent
$w_{i},w_{i+1}$ are related by a $G_{4}^{3}$-relation and all $w_{j}$ are realisable. \label{mtr}
\end{thm}

Unfortunately, we don't know whether this theorem holds for $n$-strand braids and words in $G_{n}^{3}$
for arbitrary $n\ge 4$.

Nevertheless, the construction from Section \ref{sct} works for arbitrary $n\ge 4$.

This theorem \ref{mtr} is a partial case of a phenomenon widely known in topology and group theory: by adding new
types of diagrams and moves (or by adding new generators and relations) we impose no further restrictions on
existing equivalence classes.

Other known particular cases are as follows:
\begin{enumerate}

\item Two classical link diagrams equivalent as virtual link diagrams give rise to the same link
(\cite{Kup,ClProj}.

\item Two classical braid diagrams (braid-words) equal in the virtual braid word are equal in
the classical braid group. \cite{FRR,BraidProj}.

\end{enumerate}

The key lemma for the key theorem is the following:

\begin{lm}
There exists a way of identifying ``bad letters'' in any word $w\in W$ (the other
letters are called ``good'') and the projection $pr: W\to W$ such that:

\begin{enumerate}

\item $pr$ consists of deleting all bad letters of $w$;

\item a word $w\in W$ consists of good letters only if and only if $w$ is (algebraically) realisable.

\item If $w,w'$ are related by a $G_{n}^{3}$-relation then  $pr(w),pr(w')$ are
either related by a $G_{4}^{3}$-relation or coincide.

\end{enumerate}

\end{lm}

From this lemma, one immediately gets the main theorem: we may define the {\em stable
projection}, i.e., such $pr^{st}=pr^{k}$ such that $pr^{k}(w)=pr^{k+1}(w)$. This stable
projection takes a sequence of $G_{4}^{3}$-relations between $w,w'$ into a sequence
of relations with all intermediate words realisable.

To prove the lemma, in turn, we shall use (triple) indices.

First, given four points $\{z_{1},z_{2},z_{3},z_{4}\}$ on the plane with no three collinear ones, we define
$(i,j,k)$ to be $(+1)$ if the triangle $(z_{i},z_{j},z_{k})$ is oriented clockwise
and $(-1)$ otherwise. Obviously, $(i,j,k)=(j,k,i)=(k,i,j)=-(j,i,k)=-(i,k,j)=-(k,j,i)$.

If $z_{i},z_{j},z_{k}$ are (almost) collinear in this order,
we easily see that for each $p\neq i,j,k: (i,j,p)=(i,k,p)=(j,k,p)$.

Now, the crucial lemma for proving the main theorem is as follows.
\begin{lm}
Let $w=u (a_{ijk}a_{ijl}a_{ikl}a_{jkl}) v\in W$ be a word; let
$w'= u(a_{jkl}a_{ikl}a_{ijl}a_{ijk})v$ be obtained from $w$ by applying a $G_{4}^{3}$-relation.

Then among the four letters $a_{m}$ indicated in the LHS there are either $0$ or $1$, or $4$
good ones. The number of good words in the RHS is the same; moreover, if this number is $1$, then
the good generators in the LHS and RHS are the same (say, $a_{ijl}$).

Moreover, the status of all letters in $u$ or in $v$ is the same for $w$ and for $w'$.
\label{lemm2}
\end{lm}

Before proving this (straightforward but rather boring) lemma, we formulate the lemmas for two simpler
relations in $G_{4}^{3}$, which are more or less obvious.

\begin{lm}
Let $w=u (a_{ijk}a_{ijk}) v,w'=uv$. Then either both $a_{ijk}a_{ijk}$ in $w$ are good or both are bad.
Moreover, the status of all letters in $u$ or in $v$ is the same for $w$ and for $w'$.
\label{lem2}
\end{lm}
Indeed, after applying $a_{ijk}$, none of the indices $a_{ijp},a_{ikp},a_{jkp},p\neq i,j,k$ changes and after
applying $a_{ijk}$ twice all triple indices return to the initial position.

\begin{lm}
Let $w=u a_{m} a_{m'} v, w'= u a_{m'} a_{m} v$, where $m$ and $m'$ share no more than
one common index. Then the status of the $a_{m}$ in $w$ is the same as that of $a_{m}$ in $w'$,
the same for $a_{m},a_{m'}$ in the LHS and RHS, and
Moreover, the status of all letters in $u$ or in $v$ is the same for $w$ and for $w'$.
\label{lem3}
\end{lm}
Indeed, since $m$ and $m'$ share no two common indices, the occurency of $a_{m}$ before $a_{m'}$
does not affect the fact whether $a_{m'}$ is good and vice versa.

\vspace{0.5cm}

Having the above three lemmas, we see that after deleting the bad letters from $w$ and from $w'$ we see that $f(w)$ and $f(w')$ either coincide
or differ by the $G_{4}^{3}$-relation, which suffices to prove the main theorem.

\subsection{Proof of Lemma \ref{lemm2}}

Assume that for some $t$ and some small $\eps$
the points $z_{j}(t-\eps)$ (obtained by applying the word $u$ to the initial configuration)
can undergo $a_{ijk},a_{ijl},a_{ikl},a_{jkl}$ in this order when passing to $z_{j}(t+\eps)$.

Note that the set $\{i,j,k,l\}=\{1,2,3,4\}$ and there are no other indices in consideration.

We claim that:
{\em if at least two of the above letters in the word $w$ are realisable then all four are realisable}.

We shall prove only one case of the above lemma, {\em when one of the two generators is the first one $a_{ijk}$}.
First, since we have lemmas \ref{lem2},\ref{lem3}, we can transform the main relation for $G_{n}^{3}$ by cyclic
permutations to get the desired partial case; besides, the proof in all other cases follows exactly the same procedure.

Hence, we assume $w=u (a_{ijk}a_{ijl}a_{ikl}a_{jkl}) v\in W$; let us analyse the triple indices
that we get after applying $u$ to the initial indices.

{\bf Case 1.} Assume $a_{ijk},a_{ijl}$ in $w$ are realisable.

There are three subcases:

$a)$ $a_{ijk}$ is realisable with respect to the central letter $j$.

Then after applying the subword $u$ of $w$, we get
$(i,j,l)=(i,k,l)=(j,k,l)$; without loss of generality we may think
that they are all $1$. Denote by $(i,j,k)_{0}$ the triple index $(i,j,k)$ after
applying $u$ and denote by $(i,j,k)_{1}=-(i,j,k)_{0}$ the triple index $(i,j,k)$ after applying
$u\cdot a_{ijk}$.

Now, we see that $(i,l,k)=-(i,k,l)=-1$ and $(j,l,k)=-(j,k,l)=-1$, which means that
the letter $a_{ijl}$ in $w$ can be realisable only either with respect to $i$ or with respect to $j$.

If it is realisable with respect to $j$, then we are forced to have $(i,j,k)_{1}=-1$.
After applying $a_{ikl}$, our indices become
$(i,j,k)' = (i,j,k)_{1}=-1, (i,j,l)'=-(i,j,l)=-1,(i,k,l)'=-(i,k,l)=-1,(j,k,l)'=(j,k,l)=1$.
Then one verifies straightforwardly that the next two letters $a_{ikl}$, $a_{jkl}$ are realisable
both with respect to $k$.

Also, one immediately checks that when we start with $(i,j,l)=(i,k,l)=(j,k,l)=(i,j,k)=1$,
the four letters $a_{jkl},a_{ikl},a_{ijl},a_{ijk}$ in $w'$ will be realisable with respect
to $k,k,j,j$.

Certainly, there is no problem with checking triples of indices containing some $p\notin \{i,j,k\}$.

b)
Now, assume that $a_{ijl}$ is realisable with respect to $i$. Then we may rewrite the four
letters in $w$ with respect to the ordering $l<i<j<k$:
$a_{ijk}a_{lij}a_{lik}a_{ljk}$. Here one checks straightforwardly that these four letters are
realisable with respect to $j,i,i,j$ respectively, and the same happens
for the letters $a_{ljk}a_{lik}a_{lij}a_{ijk}$ in $w'$.

c) $a_{ijk}=a_{ikj}$ is realisable with respect to the middle letter $k$. This contradicts
$(i,k,l)=(j,k,l)$.

Now, for the remaining two cases (where $a_{ijk}$ is realisable together with $a_{ijl}$ or with $a_{jkl}$
we just indicate the possible total orderings of letters meaning that the generator is realisable with
respect to the middle letter).

Here we just give the final answer without going to the detail. It suffices to observe that
any realisability of four letters is related
to some total ordering of letters $i,j,k,l$ meaning that for some ordering (say, $i<j<k<l$
each of the four words is realisable with respect to the letter staying in the middle)
and then to check that if some two of the four letters realisable then this exactly fits
into one of the ordering.

Namely, if two of the four letters of $w$ are realisable, then we get the realisability of all four letters
with respect to one of the following orderings:
$i<j<k<l$ (or $l<k<j<i$);
$j<k<l<i$ (or $i<l<k<j$);
$k<l<i<j$ (or $j<i<l<k$)
and the letters for $w'$ are realisable with respect to exactly the same ordering.

\subsection{Restoring a braid from a relaisable word}

\label{sct}

Our next goal is to describe how to restore a cylindrical braid from a word from an algebraically
realisable $w\in W$.

Consider the alphabet consisting of $(n-1)(n-2)$ letters $b_{ij}, i,j=1, \cdots, n-1,i\neq j$.

Namely, let us consider a realisable word $w$.

We shall concentrate only on generators $a_{ijn}$ containing $n$ in its index. They correspond to
trisecants passing through $z_{n}(t)$. Since the word $W$ is realisable, for each trisecant
$a_{ijn}$ we can identify which of the three letters $z_{i}(t),z_{j}(t)$ or $z_{n}(t)$ lies in the centre
(by looking at indices $(i,j,l),(i,n,l),(j,n,l)$) for $l\neq i,j,n$.

Indeed, take any $l\neq i,j,n$: if $(n,i,l)=(n,j,l)=(i,j,l)$ then $z_{i}$ is closer to
$z_{n}$ than $z_{j}$; if $(n,i,l)=(n,j,l)=(j,i,l)$ then $z_{j}$ is closer to $z_{n}$ than $z_{i}$.

We disregard letters with $n$ in the centre.

Then we get a word $pr_{n}(w)$ in the alphabet $\{b_{ij}\}$ by associating $b_{ij}$ with
those letters $a_{nij}$ in $w$ having $i$ in the centre. We call those letters $a_{njk}$ in $w$ which
give rise to $b_{jk}$ {\em surviving letters}.

This means that we can restore the picture how the points $z_{1},\cdots, z_{n-1}$ wind around $z_{n}$.

This allows one to draw a picture of a cylindrical braid.

Indeed, we may assume that $z_{n}(t)$ is constantly $0$, and consider the rays
from $0$ to $z_{j}$.

We have a cyclic ordering of angles $arg(z_{j}),j=1,\cdots, n$.

Once we pass through a trisecant $(z_{n}(t),z_{i}(t),z_{j}(t))$ where $z_{n}$ is not in the centre,
the two angles ($arg(z_{i}(t)),arg(z_{j}(t))$) switch, and we can draw the corresponding crossing on the cylinder.

Note that the two switching angles are adjacent on the circle because of the realisability
of the word $w$.

The only ambiguity in restoring the initial braid consists of the full twists.

It remains to note that if two realisable words $w,w'$ are related by a $G_{n}^{3}$-relation
then $pr_{n}(w)$ and $pr_{n}(w')$ give rise to isotopic cylindrical braids.

The most interesting case is related to the relations of length $8$, say,
$a_{nij}a_{nik}n_{njk}a_{ijk}=a_{ijk}a_{njk}a_{nik}a_{nij}$. Then it is obvious
to check that among the three letters $a_{nij},a_{nij},a_{njk}$ the number of
surviving ones is either one or three; if exactly one survives then
the word in $\{b_{ij}\}$ corresponding to this relation does not change;
if three letters survive then this gives rise to the third Reidemeister move.

Q.E.D.

\section{The induction step}

In this section we shall prove the following
\begin{thm}
If Theorem \ref{maintheorem} is true for $n=m, m\ge 4$, then it is true for $n=m+1$.
\label{induction}
\end{thm}

Obviously, this proves the main theorem.

Let us prove Theorem \ref{induction}.

We shall start with the case $m=4$, and we'll see that the general case will follow from it.

Assume $\beta$ is a $5$-strand braid and $f(\beta)=1$. Let
$\beta'$ be the braid obtained from $\beta$ by deleting the $(m+1)$-st strand.
Then we can consider the standard
projection $pr: G_{m+1}^{3}\to G_{m}^{3}$ (in our case $G_{5}^{3}\to G_{4}^{3}$) which forgets
the last strand (takes all generators $a_{i,j,m+1}$ to $1$). The image has to be trivial,
so, the braid $\beta'$ has to be homotopic to the power of full twist.

Without loss of generality, we may assume that $\beta'$ is trivial and
we may homotope $\beta$ to the position where all strands of $\beta'$ are immovable.

From now on, we assume that the first $m$ strands of $\beta$ are uniformly distributed
along the circle, i.e.,  $z_{j}(t) = exp( \frac{2\pi {j} i }{m+1}), j=1,\cdots, m$,
and the last strand $z_{m+1}(t)$ is moving somehow so that $z_{m+1}(0)=z_{m+1}(1)= 1$.

Let us consider the word $f(\beta)$. Since the first $m$ strands do not move, the word
$f(\beta)$ will contain only generators $a_{j,k,m+1}$.

By {\em principal segments} we shall mean the segments $s_{i}$ connecting points $z_{i}$ to $z_{i+1},i=1,\cdots,n-1$.
It is clear that $\C^{1} \backslash (s_{1}\cup s_{2}\cup \cdots \cup s_{n-1})$ is simply connected

Our braid gives rise to a loop $\gamma$ in the punctured plane $\C \backslash \{z_{1},\cdots, z_{m}\}$.
We assume $\gamma$ to be generic with respect to the lines connecting $z_{j},z_{k}$.
Denote the class of $\gamma$ in the corresponding homotopy group by $[\gamma]$, and denote by
$w$ the word $f(\beta)$. Assume that $\gamma$ realises the minimal length of $w$. In particular this means
that there are no consecutive intersections of $\gamma$ with the same segment.

If $\gamma$ is non-trivial and the word $w$ is non-empty then it contains an intersection
with at least one segment $s_{j}$. Hence, the word $w$ contains instances of $a_{j,j+1,n}$.
One readily checks that
the word $w$ is not trivial in $G_{n+1}^{3}$ because the letters $a_{j,j+1,n}$ can not be
cancelled with each other (see Manturov-Nikonov indices in \cite{Gnk}): for some adjacent indices
to be the same, some instance of $a_{ijk}$ with all $i,j,k\neq n$ is needed.

\section{A monomorphic map $PB_{n}\to G_{n+1}^{3}$.}

One can easily emberd $n$-stand braids into $(n+1)$-stand braids by adding the $(n+1)$-st strand
``at the infinity''.

This leads to a homomorphism $PB_{n}\to PB_{n+1}\to G_{n+1}^{3}$. It is easy to see that this homomorphism
has no kernel: one just checks that the image of full twists $tw^{m}$ is non-trivial.

In fact, an upgrade of the $G_{n}^{3}$ group where the braid group embeds was suggested in
\cite{KM}, see also \cite{book}.

However, in that paper Artin's generators were ``explicitly present''; in our approach
we deal only with $G_{n}^{3}$ with no further upgrades (orientation or coorientation, order of points
on the line etc.)

\section{Two important questions}

The proof method for $G_{n}^{4}$ is widely known in low-dimensional topology: having some ``good'' (classical,
realisable) representatives (diagrams, words), one can try to project all representatives to good ones in order
to prove that usual equivalence between good objects is the same as equivalence by means of good diagrams.

The first question is: can we find ``realisability projection'' for $PB_{n}$ for arbitrary $n$

The second question is whether one can find a faithful representation of $G_{n}^{3}$ or to find a way to use
the groups $G_{n}^{3}$ in order to get a faithful representation for the Artin braid group.

\end{document}